\documentclass[11pt]{article}
\usepackage{epic,latexsym,amssymb}
\usepackage{tikz}

\usepackage{tikz,epsf}

\usepackage{amsfonts}
\usepackage{amscd}
\usepackage{amsmath}
\usepackage{graphicx}
\usepackage{color}
\usepackage{caption,subcaption}

\newtheorem{thm}{Theorem}
\newtheorem{ob}[thm]{Observation}

\newtheorem{lem}[thm]{Lemma}

\newtheorem{cor}[thm]{Corollary}

\newcommand{\diam}{{\rm diam}}

\newcommand{\sta}{\mathop{\mathrm{sta}}}

\newcommand{\cP}{\mathcal{P}}

\newcommand{\cT}{\mathcal{T}}
\newcommand{\cO}{\mathcal{O}}

\newcommand{\proof}{\noindent\textbf{Proof. }}

\newcommand{\qed}{$\Box$}

\newcommand{\1}{\vspace{0.1cm}}

\newcommand{\gt}{\gamma_t}

\newcommand{\gnt}{\gamma_{\rm nt}}

\newcommand{\vertex}{\node[vertex]}
\tikzstyle{vertex}=[circle, draw, inner sep=0pt, minimum size=6pt]

\newcommand{\QEDmark}{\mbox{\textsc{qed}}}
\newcommand{\proofStarter}[1]{\textsc{#1} }

\newcommand{\pc}{{\rm pc}}

\def\vertex(#1){\put(#1){\circle*{2}}}
\def\vertexo(#1){\put(#1){\circle{2}}}
\def\vert(#1){\put(#1){\circle*{1.5}}}
\def\verto(#1){\put(#1){\circle{1.5}}}
\def\lab(#1)#2{\put(#1){\makebox(0,0)[c]{#2}}}
\begin{document}

\title{Matchings and Path Covers with  applications to \\ Domination in Graphs}
\author{Michael A. Henning\thanks{Research supported in part by the South African National Research Foundation and the University of Johannesburg} \, and Kirsti Wash\thanks{Research supported in part by the University of Johannesburg}\\
\\
Department of Pure and Applied Mathematics\\
University of Johannesburg \\
Auckland Park, 2006, South Africa \\
Email: mahenning@uj.ac.za\\
Email: kirsti.wash@trincoll.edu
}

\date{}
\maketitle

\begin{abstract}
Let $G$ be a graph with no isolated vertex. A matching in $G$ is a set of edges that are pairwise not adjacent in $G$, while the matching number, $\alpha'(G)$, of $G$ is the maximum size of a matching in $G$. The path covering number, $\pc(G)$, of $G$ is the minimum number of vertex disjoint paths such that every vertex belongs to a path in the cover. We show that if $G$ has order~$n$, then $\alpha'(G) + \frac{1}{2}\pc(G) \ge \frac{n}{2}$ and we provide a constructive characterization of the graphs achieving equality in this bound.
It is known that $\gamma(G) \le \alpha'(G)$ and $\gamma_t(G) \le \alpha'(G) + \pc(G)$, where $\gamma(G)$ and $\gamma_t(G)$ denote the domination and the total domination number of $G$. As an application of our result on the matching and path cover numbers, we show that if $G$ is a graph with $\delta(G) \ge 3$, then $\gamma_t(G) \le \alpha'(G) + \frac{1}{2}(\pc(G) - 1)$, and this bound is tight.
A set $S$ of vertices in $G$ is a neighborhood total dominating set of $G$ if it is a dominating set of $G$ with the property that the subgraph induced by the open neighborhood of the set $S$ has no isolated vertex. The neighborhood total domination number, $\gnt(G)$, is the minimum cardinality of a neighborhood total dominating set of $G$. We observe that $\gamma(G) \le \gnt(G) \le \gamma_t(G)$. As a further application of our result on the matching and path cover numbers, we show that if $G$ is a connected graph on at least six vertices, then $\gnt(G) \le \alpha'(G) + \frac{1}{2}\pc(G)$ and this bound is tight. 
%
\end{abstract}


\newpage
\section{Introduction}

Two distinct edges in a graph $G$ are \emph{independent} if they are not adjacent in $G$. A \emph{matching} in $G$ is a set of (pairwise) independent edges, while a matching of maximum cardinality is a \emph{maximum matching}. The number of edges in a maximum matching of a graph $G$ is called the \emph{matching number} of $G$, denoted by $\alpha'(G)$. Let $M$ be a specified matching in a graph $G$. A vertex $v$ of $G$ is an \emph{$M$-matched vertex} if $v$ is incident with an edge of $M$; otherwise, $v$ is an \emph{$M$-unmatched vertex}.
%
An $M$-\emph{alternating path} of $G$ is a path whose edges are alternately in $M$ and not in $M$.
An \emph{$M$-augmenting path} is an $M$-alternating path that begins and ends with $M$-unmatched vertices.
Given an $M$-augmenting path $P$, we denote by $M \triangle E(P)$ the symmetric difference of $M$ and $E(P)$, i.e.,  $M \triangle E(P) = (M \setminus E(P)) \cup (E(P) \setminus M)$.
Matchings in graphs are extensively studied in the literature (see, for example, the excellent survey articles by Plummer~\cite{Pl03} and Pulleyblank~\cite{Pu95}).

A \emph{path cover} in a graph $G$ is a collection of vertex disjoint paths such that every vertex belongs to exactly one path. The cardinality of a minimum path cover is the \emph{path covering number} of $G$ which we denote by $\pc(G)$.

A \emph{dominating set} in a graph $G$ is a set $S$ of vertices of $G$ such that every vertex in $V(G) \setminus S$ is adjacent to at least one vertex in $S$. The \emph{domination number} of $G$, denoted by $\gamma(G)$, is the minimum cardinality of a dominating set of $G$. A \emph{total dominating set}, abbreviated a TD-set, of a graph $G$ with no isolated vertex is a set $S$ of vertices of $G$ such that every vertex in $V(G)$ is adjacent to at least one vertex in $S$. The \emph{total domination number} of $G$, denoted by $\gt(G)$, is the minimum cardinality of a TD-set of $G$. The literature on the subject of domination parameters in graphs up to the year 1997 has been surveyed and detailed in the two books~\cite{hhs1, hhs2}.
Total domination is now well studied in graph theory. For a recent book on the topic, see~\cite{HeYe_book}. A survey of total domination in graphs can also be found in~\cite{He09}.

Arumugam and Sivagnanam~\cite{ArSi11} introduced and studied the concept of neighborhood total domination in graphs. A \emph{neighbor} of a vertex $v$ is a vertex different from $v$ that is adjacent to $v$. The \emph{neighborhood of a set} $S$ is the set of all neighbors of vertices in $S$. A \emph{neighborhood total dominating set}, abbreviated NTD-set, in a graph $G$ is a dominating set $S$ in $G$ with the property that the subgraph induced by the open neighborhood of the set $S$ has no
isolated vertex. The \emph{neighborhood total domination number} of $G$, denoted by $\gnt(G)$, is the minimum cardinality of a NTD-set of $G$.
%
Since every TD-set is a NTD-set, and since every NTD-set is a dominating set, we have the following observation first observed by Arumugam and Sivagnanam in~\cite{ArSi11}.

\begin{ob}{\rm (\cite{ArSi11})}
If $G$ is a graph with no isolated vertex, then $\gamma(G) \le \gnt(G) \le \gt(G)$.
 \label{relate1}
\end{ob}

By Observation~\ref{relate1}, the neighborhood total domination number is squeezed between arguably the two most important domination parameters, namely the domination number and the total domination number.
%

\subsection{Terminology and Notation}

For notation and graph theory terminology not defined herein, we refer the reader to~\cite{hhs1}. Let $G$ be a graph with vertex set $V(G)$ of order~$n = |V(G)|$ and edge set $E(G)$ of size~$m = |E(G)|$, and let $v$ be a vertex in $V$. We denote the \emph{degree} of $v$ in $G$ by $d_G(v)$. The minimum degree among the vertices of $G$ is denoted by $\delta(G)$.
A vertex of degree one is called a \emph{leaf} and its neighbor a \emph{support vertex}.
%
%
For a set $S \subseteq V$, the subgraph induced by $S$ is denoted by $G[S]$.

A \emph{cycle} and \emph{path} on $n$ vertices are denoted by $C_n$ and $P_n$, respectively. A \emph{star} on $n \ge 2$ vertices is a tree with a vertex of degree~$n-1$ and is denoted by $K_{1,n-1}$. A \emph{double star} is a tree containing  exactly two vertices that are not leaves (which are necessarily adjacent). A \emph{subdivided star} is a graph obtained from a star on at least two vertices by subdividing each edge exactly once. We note that the smallest two subdivided stars are the paths $P_3$ and $P_5$.

The \emph{open neighborhood} of $v$ is the set $N_G(v) = \{u \in V \, | \, uv \in E\}$ and the \emph{closed neighborhood of $v$} is $N_G[v] = \{v\} \cup N_G(v)$. For a set $S\subseteq V$, its \emph{open neighborhood} is the set $N_G(S) = \bigcup_{v \in S} N_G(v)$, and its \emph{closed neighborhood} is the set $N_G[S] = N_G(S) \cup S$. If the graph $G$ is clear from the context, we simply write $\text{d}(v)$, $N(v)$, $N[v]$, $N(S)$ and $N[S]$ rather than $d_G(v)$, $N_G(v)$, $N_G[v]$, $N_G(S)$ and $N_G[S]$, respectively. As observed in~\cite{HeRa13} a NTD-set in $G$ is a set $S$ of vertices such that $N[S] = V$ and $G[N(S)]$ contains no isolated vertex.

A \emph{rooted tree} distinguishes one vertex $r$ called the \emph{root}. For each vertex $v \ne r$ of $T$, the \emph{parent} of $v$ is the neighbor of $v$ on the unique $(r,v)$-path, while a \emph{child} of $v$ is any other neighbor of $v$. A \emph{descendant} of $v$ is a vertex $u$ such that the unique $(r,u)$-path contains $v$. Let $C(v)$ and $D(v)$ denote the set of children and descendants, respectively, of $v$, and let $D[v] = D(v) \cup \{v\}$. A \emph{non}-\emph{leaf} of a tree $T$ is a vertex of $T$ of degree at least~$2$ in $T$.

\subsection{Known Results}

Bounds relating the domination number and the matching number are studied, for example, in~\cite{BoCo79,CoHeSl79}. As a consequence of a result due to Bollob\'{a}s and Cockayne~\cite{BoCo79}, the
domination number of every graph with no isolated vertex is bounded above by its matching number.

\begin{thm}{\rm (\cite{BoCo79})}
For every graph $G$ with no isolated vertex, $\gamma(G) \le
\alpha'(G)$. \label{extnbh}
\end{thm}

The total domination number versus the matching number in a graph has been studied in several papers (see, for example,~\cite{DeLa07,HeKaShYe08,HeYe06,HeYe11,HeYe13,HeYe_chapter,OWest10,ShChCh10,WaKaSh09} and elsewhere). Unlike the domination number, the total domination number and the matching number of a graph are generally incomparable, even for arbitrarily large, but fixed (with respect to the order of the graph), minimum degree as shown in~\cite{HeKaShYe08}. The following upper bound on the total domination in terms of its matching number and path covering number is presented in~\cite{DeLa07}.

\begin{thm}{\rm (\cite{DeLa07})}
For every graph $G$ with no isolated vertex, $\gamma_t(G) \le \alpha'(G) + \pc(G)$, and this bound is tight.
 \label{Tdom_pc}
\end{thm}


The following upper bound on the neighborhood total domination number of a connected graph in terms of its order is established in~\cite{HeRa13}.

\begin{thm}{\rm (\cite{HeRa13})}
Let $G$ be a connected graph of order~$n \ge 3$. Then, $\gnt(G) \le (n+1)/2$ with equality if and only if $G = C_5$ or $G$ is a subdivided star.
 \label{t:oddn}
\end{thm}

\subsection{The Family $\cT$}

By a \emph{weak partition} of a set we mean a partition of the set in which some of the subsets may be empty. For our purposes we define a \emph{labeling} of a tree $T$ as a weak partition $S = (S_A,S_B,S_C)$ of $V(T)$. We will refer to the pair $(T,S)$ as a \emph{labeled tree}. The \emph{label} or \emph{status} of a vertex
$v$, denoted $\sta(v)$, is the letter $x \in \{A,B,C\}$ such that $v \in S_x$. We now define the following family of trees.

Let $\cT$ be the family of labeled trees that: (i) contains $(P_3,S_0^*)$ where $S_0^*$ is the labeling that assigns status $A$ to the central vertex of $P_3$ and status~$B$ to the two leaves; and (ii) is closed under the four operations $\cO_1$, $\cO_2$, $\cO_3$ and $\cO_4$ that are listed below, which extend the tree $T'$ to a tree $T$ by attaching a tree to the vertex $v \in V(T')$, called the \emph{attacher} of~$T'$.

\begin{itemize}
\item \textbf{Operation $\cO_1$.} Let $v$ be a vertex with $\sta(v) = A$. Add a vertex $u_1$ and the edge $vu_1$. Let $\sta(u_1)=B$.

\item \textbf{Operation $\cO_2$.} Let $v$ be a vertex with $\sta(v) = B$. Add a path $u_1u_2$ and the edge $vu_1$. Change the status of $v$ from $B$ to $C$, and so $\sta(v) = C$, and let $\sta(u_1) = A$ and $\sta(u_2) = B$.

\item \textbf{Operation $\cO_3$.} Let $v$ be an arbitrary vertex of $T'$. Add a path $u_1u_2u_3$ and the edge $vu_2$. Let    $\sta(u_1) = \sta(u_3) = B$ and $\sta(u_2) = A$. Further, if $\sta(v) = B$, then change the status of $v$ from $B$ to $C$.

\item \textbf{Operation $\cO_4$.} Let $v$ be a vertex with $\sta(v) = A$. Add a path $u_1u_2u_3u_4u_5$ and the edge $vu_3$. Let    $\sta(u_1) = \sta(u_5) = B$, $\sta(u_2) = \sta(u_4) = A$, and $\sta(u_3) = C$.

\end{itemize}

The four operations $\cO_1$, $\cO_2$, $\cO_3$  and $\cO_4$ are illustrated in Figure~\ref{f:O1O2}.

\setlength{\unitlength}{0.9mm}
\begin{figure}[htb]
\begin{center} \itshape

\begin{picture}(80,20)(-42,-12)
\put(-32,0){$\cO_1$:} \vertex(-12,0) \lab(-12,3){A} \vertex(0,0) \lab(0,3){B}
\put(-12,0){\line(1,0){12}}
\put(-12,0){\circle{22}}
\end{picture}

\begin{picture}(80,20)(-42,-12)
\put(-32,0){$\cO_2$:} \vertex(-12,0) \lab(-12,3){B} \vertex(0,0) \lab(0,3){A} \lab(12,3){B} \vertex(12,0)
\put(-12,0){\line(1,0){24}} \put(-12,0){\circle{22}}
\end{picture}

\begin{picture}(80,20)(-42,-12)
\put(-32,0){$\cO_3$:} \vertex(-12,0) \vertex(0,0) \vertex(12,6) \vertex(12,-6)
\lab(-12,3){{\small any}}
\lab(0,3){A} \lab(15,6){B} \lab(15,-6){B}
\put(-12,0){\line(1,0){12}}
\put(0,0){\line(2,1){12}}
\put(0,0){\line(2,-1){12}}
\put(-12,0){\circle{22}}
\end{picture}

\begin{picture}(80,25)(-42,-12)
\put(-32,0){$\cO_4$:} \vertex(-12,0) \vertex(0,0) \vertex(12,6) \vertex(12,-6) \vertex(24,6) \vertex(24,-6)
\lab(-12,3){A}
\lab(0,3){C} \lab(12,9){A} \lab(12,-9){A} \lab(24,9){B} \lab(24,-9){B}
\put(-12,0){\line(1,0){12}}
\put(0,0){\line(2,1){12}}
\put(0,0){\line(2,-1){12}}
\put(12,6){\line(1,0){12}}
\put(12,-6){\line(1,0){12}}
\put(-12,0){\circle{22}}
\end{picture}

\end{center}
\vskip -0.25 cm
\caption{The four operations $\cO_1$, $\cO_2$, $\cO_3$  and $\cO_4$.}
\label{f:O1O2}
\end{figure}

We shall need the following properties of labeled trees $(T,S)$ in the family~$\cT$ that follow immediately from the way in which trees in the family~$\cT$ are constructed.

\begin{ob}
If $(T,S) \in \cT$ is a labeled tree for some labeling $S$, then $T$ has the following properties:
\\
\indent {\rm (a)} Every support vertex of $T$ belongs to $S_A$.
\\
\indent {\rm (b)} The set $S_B$ is the set of leaves of $T$.
\\
\indent {\rm (c)} Every vertex in $S_A$ has at least two neighbors in $S_B \cup S_C$.
\\
\indent {\rm (d)} If $v \in S_C$, then $N(v) \subseteq S_A$.
\label{ob:cT}
\end{ob}

\section{Main Results}

We have two immediate aims in this paper. Our first aim is to establish a lower bound relating the matching number and the path covering number of a graph in terms of the order of the graph. Our second aim is to
establish an upper bound on the neighborhood total domination number of a graph in terms of its matching and path covering numbers. More precisely, we shall prove the following results, proofs of which are presented in Sections~\ref{S:main2} and~\ref{S:applic}.

\begin{thm}
Let $T$ be a tree of order~$n \ge 3$. Then, $\alpha'(T) + \frac{1}{2}\pc(T) \ge \frac{n}{2}$ with equality if and only if $(T,S) \in \cT$ for some labeling $S$.
 \label{t:tree}
\end{thm}

\begin{thm}
If $G$ is a graph of order~$n$, then $\alpha'(G) + \frac{1}{2}\pc(G) \ge \frac{n}{2}$. Further for $n \ge 3$, $\alpha'(G) + \frac{1}{2}\pc(G) = \frac{n}{2}$ if and only if $G$ has a spanning tree $T$ such that \1
\\
\indent {\rm (a)} $(T,S) \in \cT$ for some labeling $S$.
\\
\indent {\rm (b)} $\alpha'(G) = \alpha'(T)$. \\
\indent {\rm (c)}  $\pc(G) = \pc(T)$.
 \label{t:thm2}
\end{thm}

As an application of Theorem~\ref{t:thm2}, we show that the bound of Theorem~\ref{Tdom_pc} can be improved considerably if we restrict the minimum degree of the graph to be at least three. 

\begin{cor}
If $G$ is a graph with $\delta(G) \ge 3$, then $\gamma_t(G) \le \alpha'(G) + \frac{1}{2}(\pc(G) - 1)$, and this bound is tight.
 \label{t:thm4}
\end{cor}

As a further application of Theorem~\ref{t:thm2}, we establish the following upper bound on the neighborhood total domination number of a graph in terms of its matching and path covering numbers.

\begin{cor}
If $G$ is a connected graph of order at least~$3$, then
\[
\gnt(G) \le \alpha'(G) + \frac{1}{2}\pc(G)
\]
unless $G \in \{P_3,P_5,C_5\}$ in which case $\gnt(G) = \alpha'(G) + \frac{1}{2}(\pc(G)+1)$.
 \label{t:thm3}
\end{cor}

\section{Proof of Theorem~\ref{t:tree} and Theorem~\ref{t:thm2}}
\label{S:main2}

In this section, we establish a lower bound for $\alpha'(G) + \frac{1}{2}\pc(G)$. For this purpose, we shall need the following results on matchings and path covers in trees.

\begin{lem}
If $T$ is a tree of order~$n \ge 1$, then $\alpha'(T) + \frac{1}{2}\pc(T) \ge \frac{n}{2}$.
 \label{l:tree1}
\end{lem}
\proof We proceed by induction on the order $n \ge 1$ of a tree $T$. If $n \in \{1,2,3\}$, then $T = P_n$ and the result follows readily. This establishes the base case. Let $n \ge 4$ and assume that if $T'$ is a tree of order $n'$ where $n' < n$, then $\alpha'(T') + \frac{1}{2}\pc(T') \ge \frac{n'}{2}$. Let $T$ be a tree of order $n$.
If $T$ is a star, then $\alpha'(T) = 1$ and $\pc(T) = n-2$, and so $\alpha'(T) + \frac{1}{2}\pc(T) = \frac{n}{2}$.
Hence we may assume that $\diam(T) \ge 3$, for otherwise the desired result follows.
Let $P$ be a longest path in $T$ and suppose that $P$ is an $(r,u)$-path. Necessarily, $r$ and $u$ are leaves in $T$. We now root the tree $T$ at the vertex~$r$. Let $v$ be the parent of $u$, and let $w$ be the parent of $v$ in the rooted tree $T$.

Let $\cP$ be a minimum path cover in $T$ and let $P_v$ be the path in $\cP$ containing~$v$. Thus, $|\cP| = \pc(T)$. By the minimality of the path cover $\cP$ in $T$, we may assume, renaming the children of $v$ if necessary, that $u \in V(P_v)$.

Suppose that
%
$d_T(v) = 2$. In this case, let $T' = T - \{u,v\}$ and let $T'$ have order $n'$, and so $n' = n-2 \ge 2$. Applying the inductive hypothesis to $T'$, we note that $\alpha'(T') + \frac{1}{2}\pc(T') \ge \frac{n'}{2}$. Every matching in $T'$ can be extended to a matching in $T$ by adding to it the edge $uv$, implying that $\alpha'(T) \ge \alpha'(T') + 1$. By our earlier assumption, the vertex $u \in V(P_v)$. If $w \in V(P_v)$, then replacing the path $P_v$ in $\cP$ with the path $P_v - \{u,v\}$ produces a path cover in $T'$ of size~$|\cP|$. If $w \notin V(P_v)$, then removing the path $P_v$ from $\cP$ produces a path cover in $T'$ of size~$|\cP| - 1$. In both cases, we produce a path cover in $T'$ of size at most~$|\cP|$, implying that $\pc(T') \le \pc(T)$. Therefore, $\alpha'(T) + \frac{1}{2}\pc(T) \ge \alpha'(T') + \frac{1}{2}\pc(T') + 1 \ge \frac{n'}{2} + 1 = \frac{n}{2}$.
%
%
Hence we may assume that $d_T(v) \ge 3$, for otherwise the desired result holds.

Let $C(v) = \{u_1,u_2,\ldots,u_k\}$ denote the children of $v$, where $u = u_1$. By assumption, $k \ge 2$. We now let $T'$ be the tree obtained from $T$ by deleting $v$ and all children of $v$; that is, $T' = T - D[v]$. Let $T'$ have order $n'$, and so $n' = n - k - 1$. Applying the inductive hypothesis to $T'$, we note that $\alpha'(T') + \frac{1}{2}\pc(T') \ge \frac{n'}{2}$. Every matching in $T'$ can be extended to a matching in $T$ by adding to it the edge $uv$, implying that $\alpha'(T) \ge \alpha'(T') + 1$. By the minimality of the path cover $\cP$ in $T$, the path $P_v$ in $\cP$ contains either the vertex~$w$ or a child of $v$ different from $u_1$. If $w \in V(P_v)$, then we can remove from $\cP$ both the path $P_v$ and the path in $\cP$ that consists only of the vertex $u_2$ and replace these two paths with the path $P_v - \{u,v\}$ and the path $u_1vu_2$. Hence, we can choose $\cP$ so that the path $P_v$ in $\cP$ that contains $v$ is the path $u_1vu_2$. We note that if $k \ge 3$, then each child $u_i$ of $v$ belongs to a path in $\cP$ that consists only of the vertex $u_i$. Removing the $k-1$ paths from $\cP$ that contain vertices in $D[v]$, we produce a path cover in $T'$ of size $|\cP| - (k-1) = \pc(T) - k + 1$, implying that $\pc(T') \le \pc(T) - k + 1$. Hence,

\[
\begin{array}{lcl}
\alpha'(T) + \frac{1}{2}\pc(T) & \ge & (\alpha'(T') + 1) + \frac{1}{2}(\pc(T') + k - 1) \1  \\
& = & (\alpha'(T') + \frac{1}{2}\pc(T')) + \frac{1}{2}(k + 1)  \1 \\
& \ge & \frac{n'}{2} + \frac{k+1}{2}  \1 \\
& = & \frac{n - k - 1}{2} + \frac{k+1}{2} \1
\\
& = & \frac{n}{2}.
\end{array}
\]

\noindent
This completes the proof of Lemma~\ref{l:tree1}.~\qed

\medskip
Let $(T,S) \in \cT$ be a labeled tree for some labeling $S$. Then there is a sequence of labeled trees $(T_0,S_0), (T_1,S_1), \ldots, (T_k,S_k)$ such that $(T_0,S_0) = (P_3,S_0^*)$, $(T_k,S_k) = (T,S)$, and if $k \ge 1$, then for $i \in \{1,\ldots,k\}$, the labeled tree $(T_i,S_i)$ can be obtained from $(T_{i-1},S_{i-1})$ by one of the operations $\cO_1$, $\cO_2$, $\cO_3$ or $\cO_4$. We call the number of terms in such a sequence of labeled trees that is used to construct $(T,S)$, the \emph{length} of the sequence. In particular, the above sequence has length~$k$. We remark that a sequence of labeled trees used to construct $(T,S)$ is not necessarily unique. Further, the length of such sequences may differ. We shall need the following properties of trees in the family~$\cT$.

\begin{lem}
If $(T,S) \in \cT$ is a labeled tree for some labeling $S$, then $T$ has the following properties.
\\
\indent {\rm (a)}  If $w \in S_A$ and $M$ is an arbitrary maximum matching in $T$, then $w$ is $M$-matched. \\
\indent {\rm (b)} If $w \in S_B \cup S_C$, then there exists a maximum matching $M$ in $T$ such that $w$ is \\ \hspace*{0.9cm} $M$-unmatched. \\
\indent {\rm (c)} $\alpha'(T) + \frac{1}{2}\pc(T) = \frac{n}{2}$.
 \label{lem:cT}
\end{lem}
\proof Let $(T,S) \in \cT$ be a labeled tree for some labeling $S$. We proceed by induction on the length, $k \ge 0$, of a sequence used to construct $(T,S)$. If $k = 0$, then $T = P_3$ and $S = S_0^*$, and Properties~(a), (b) and~(c) hold. This establishes the base case. Let $k \ge 1$ and assume that if the length of the sequence used to construct a labeled tree $(T',S') \in \cT$ is less than~$k$, then Properties~(a), (b) and~(c) hold for the labeled tree $(T',S')$. Let $(T,S) \in \cT$ and let $(T_0,S_0), (T_1,S_1), \ldots, (T_k,S_k)$ be a sequence of length~$k$ used to construct~$(T,S)$, where $(T_0,S_0) = (P_3,S_0^*)$ and $(T_k,S_k) = (T,S)$. Let $T' = T_{k-1}$ and let $S' = S_{k-1}$. Then, $(T',S') \in \cT$. Let $M$ be a maximum matching in $T$ and let $M'$ be the restriction of $M$ to the tree $T'$; that is, $M' = M \cap E(T')$. We consider four cases, depending on the operation applied to $(T',S')$ in order to obtain $(T,S)$.

\medskip
\emph{Case~1. $(T,S)$ can be obtained from $(T',S')$ by operation $\cO_1$.} Let $v$ be the attacher in $T'$ and let $u_1$ be the vertex added to $T'$ to obtain $T$. Then, $\sta(v) = A$ and $\sta(u_1) = B$. We show that $\alpha'(T') = \alpha'(T)$. Since $M$ is a maximum matching in $T$, the vertex $v$ is $M$-matched. If $vu_1 \notin M$, then $M' = M$, implying that $\alpha'(T') \ge |M'| = |M| = \alpha'(T)$. If $vu_1 \in M$, then $M' = M \setminus \{vu_1\}$ and $|M'| = |M| - 1$. Applying the inductive hypothesis to the labeled tree $(T',S') \in \cT$, the matching $M'$ is not a maximum matching in $T'$ since the vertex $v \in S_A'$ is $M'$-unmatched. Hence in this case, $\alpha'(T') \ge |M'| + 1 = |M| = \alpha'(T)$. In both cases, $\alpha'(T') \ge \alpha'(T)$. Conversely, every matching in $T'$ is a matching in $T$, and so $\alpha'(T) \ge \alpha'(T')$. Consequently, $\alpha'(T') = \alpha'(T)$.

We show firstly that Property~(a) holds. Let $w \in S_A$. Then, $w \in V(T')$. Let $M$ be an arbitrary maximum matching in $T$ and let $M'$ be the restriction of $M$ to the tree $T'$. Suppose, to the contrary, that $w$ is $M$-unmatched. If $M = M'$, then applying the inductive hypothesis to the labeled tree $(T',S') \in \cT$ the vertex $w$ is $M'$-matched in $T'$ and therefore $M$-matched in $T$, a contradiction. Hence, $M \ne M'$, implying that $vu_1 \in M$. Thus, $M' = M \setminus \{vu_1\}$ and $|M'| = |M| - 1$. We now consider the matching $M'$. Let $P_v$ be a longest path in $T'$ starting at $v$ and whose edges are alternately not in $M'$ and in $M'$ and whose vertices are alternately in $S_A'$ and in $S_C'$. Let $P_v$ be a $(v,z)$-path. By Observation~\ref{ob:cT}, every vertex in the set $S_A'$ has at least two neighbors in $S_B' \cup S_C'$.

If $v = z$, then in the tree $T'$ the vertex $v$ has no neighbor in $S_C'$ and therefore at least two neighbors in $S_B'$. Let $u'$ be a neighbor of $v$ in $T'$ that belongs to $S_B'$. By Observation~\ref{ob:cT}(b), the vertex $u'$ is a leaf. Since $vu_1 \in M$, we note that $u'$ is $M$-unmatched. We now consider the matching $M^* = M' \cup \{u'v\}$. Since $M^*$ has size~$|M'| + 1 = |M| = \alpha'(T) = \alpha'(T')$, the matching $M^*$ is a maximum matching in $T'$. Since $w \notin \{u',v\}$ and $w$ is $M$-unmatched, the vertex $w$ is $M^*$-unmatched, contradicting Property~(a). Hence, $v \ne z$.

If $z \in S_C'$, then by the maximality of the path $P_v$ the vertex $z$ is $M'$-unmatched. Thus, $P_v$ is an $M'$-augmenting $(v,z)$-path in $T'$. We now consider the matching $M^* = M \triangle E(P_v)$. Since $M^*$ has size~$|M'| + 1 = |M| = \alpha'(T) = \alpha'(T')$, the matching $M^*$ is a maximum matching in $T'$. Since $w \notin V(P_v)$ and $w$ is $M$-unmatched, the vertex $w$ is $M^*$-unmatched, contradicting Property~(a). Hence, $z \in S_A'$.

Since $z \in S_A'$, the vertex $z$ has at least two neighbors in $S_B' \cup S_C'$. By the maximality of the path $P_v$ the vertex $z$ has no neighbor in $S_C'$ except for the vertex immediately preceding it on the path $P_v$ (that is $M$-matched to $z$). Thus, $z$ has a neighbor, $z'$ say, that belongs to $S_B'$. By Observation~\ref{ob:cT}(b), the vertex $z'$ is a leaf. We now extend the path $P_v$ by adding to it the vertex $z'$ and the edge $zz'$ to produce a new path $P$ which is an $M'$-augmenting $(v,z')$-path in $T'$. We now consider the matching $M^* = M \triangle E(P)$. Since $M^*$ has size~$|M'| + 1 = |M| = \alpha'(T) = \alpha'(T')$, the matching $M^*$ is a maximum matching in $T'$. Since $w \notin V(P)$ and $w$ is $M$-unmatched, the vertex $w$ is $M^*$-unmatched, contradicting Property~(a). Therefore, $w$ is $M$-matched. Since $w$ was chosen to be an arbitrary vertex in $S_A$, this proves that Property~(a) holds.

We show secondly that Property~(b) holds. Let $w \in S_B \cup S_C$. Suppose that $w = u_1$. Let $M^*$ be an arbitrary maximum matching in $T'$. Since $\alpha'(T') = \alpha'(T)$, the matching $M^*$ is a maximum matching in $T$ such that $w$ is $M^*$-unmatched, as desired. Hence we may assume that $w \ne u_1$, for otherwise the desired result holds. Thus, $w \in V(T')$. Applying the inductive hypothesis to the labeled tree $(T',S')$, there exists a maximum matching $M_w$ in $T'$ such that $w$ is $M_w$-unmatched. Since $\alpha'(T') = \alpha'(T)$, the matching $M_w$ is a maximum matching in $T$ such that $w$ is $M_w$-unmatched. This establishes Property~(b).

We prove next that Property~(c) holds. As observed earlier, $\alpha'(T') = \alpha'(T)$. Either $\pc(T') = \pc(T) - 1$ or $\pc(T') = \pc(T)$, implying by Lemma~\ref{l:tree1} that $\frac{n}{2} \le \alpha'(T) + \frac{1}{2}\pc(T) \le \alpha'(T') + \frac{1}{2}(\pc(T') + 1) = \frac{n'}{2} + \frac{1}{2} = \frac{n}{2}$. Hence we must have equality throughout this inequality chain. In particular, $\alpha'(T) + \frac{1}{2}\pc(T) = \frac{n}{2}$. This establishes Property~(c).

\medskip
\emph{Case~2. $(T,S)$ can be obtained from $(T',S')$ by operation $\cO_2$.} Let $v$ be the attacher in $T'$ and let $u_1u_2$ be the path added to $T'$ and $vu_1$ the added edge. Then, $\sta_{T'}(v) = B$ but $\sta_{T}(v) = C$, while $\sta(u_1) = A$ and $\sta(u_2) = B$. If $u_1u_2 \notin M$, then by the maximality of $M$, the edge $u_1v \in M$ and we can simply replace the edge $u_1v$ in $M$ with the edge $u_1u_2$. Hence we may choose $M$ so that $u_1u_2 \in M$. Thus, $M' = M \setminus \{u_1u_2\}$ is a matching in $T'$, implying that $\alpha'(T') \ge |M'| = |M| - 1 = \alpha'(T) - 1$. Every matching in $T'$ can be extended to a matching in $T$ by adding to it the edge $u_1u_2$, and so $\alpha'(T) \ge \alpha'(T') + 1$. Consequently, $\alpha'(T') = \alpha'(T) - 1$.

We show firstly that Property~(a) holds. Let $w \in S_A$. Let $M$ be an arbitrary maximum matching in $T$ and let $M'$ be the restriction of $M$ to the tree $T'$. Suppose, to the contrary, that $w$ is $M$-unmatched. If $w = u_1$, then by the maximality of $M$, the vertex $w$ is $M$-matched, a contradiction. Hence, $w \in V(T')$. Suppose $vu_1 \in M$. Then, $M' = M \setminus \{vu_1\}$. Since $\alpha'(T') = \alpha'(T) - 1$ and $M'$ has size~$|M| - 1 = \alpha'(T) - 1 = \alpha'(T')$, the matching $M'$ is a maximum matching in $T'$. Applying the inductive hypothesis to the labeled tree $(T',S') \in \cT$ the vertex $w$ is $M'$-matched in $T'$ and therefore $M$-matched in $T$, a contradiction. Therefore, $vu_1 \notin M$, implying that $u_1u_2 \in M$ and $M' = M \setminus \{u_1u_2\}$. Once again, $M'$ has size~$\alpha'(T) - 1$ and is a maximum matching in $T'$. Applying the inductive hypothesis to the labeled tree $(T',S') \in \cT$, the vertex $w$ is $M'$-matched in $T'$ and therefore $M$-matched in $T$, a contradiction. This establishes Property~(a).

We show secondly that Property~(b) holds. Let $w \in S_B \cup S_C$. Suppose that $w = u_2$. Applying the inductive hypothesis to the labeled tree $(T',S')$, there is a maximum matching $M^*$ of $T'$ such that the vertex $v$ which has status~$B$ in $(T',S')$ is $M^*$-unmatched. Thus the matching $M_w = M^* \cup \{vu_1\}$ is a maximum matching in $T$ such that $w$ is $M_w$-unmatched. Hence we may assume that $w \in V(T')$, for otherwise the desired result holds. (Possibly, $w = v$, in which case $w$ has status $B$ in $(T',S')$ and status~$C$ in $(T,S)$.)  Applying the inductive hypothesis to the labeled tree $(T',S') \in \cT$, there is a maximum matching $M^*$ of $T'$ such that the vertex $w$ is $M^*$-unmatched. Thus the matching $M_w = M^* \cup \{u_1u_2\}$ is a maximum matching in $T$ such that $w$ is $M_w$-unmatched. This establishes Property~(b).

We prove next that Property~(c) holds. As observed earlier, $\alpha'(T') = \alpha'(T) - 1$. Since the attacher vertex $v$ is a leaf in $T'$, we note that $\pc(T) = \pc(T')$. Therefore, $\alpha'(T) + \frac{1}{2}\pc(T) = (\alpha'(T') + 1) + \frac{1}{2}\pc(T') = \frac{n'}{2} + 1 = \frac{n}{2}$. This establishes Property~(c).

\medskip
\emph{Case~3. $(T,S)$ can be obtained from $(T',S')$ by operation $\cO_3$.} Let $v$ be the attacher in $T'$ and let $u_1u_2u_3$ be the path added to $T'$ and $vu_2$ the added edge. Then, $\sta(u_1) = \sta(u_3) = B$ and $\sta(u_2) = A$. Further, $\alpha'(T') = \alpha'(T) - 1$.

We show firstly that Property~(a) holds. Let $w \in S_A$. Let $M$ be an arbitrary maximum matching in $T$ and let $M'$ be the restriction of $M$ to the tree $T'$. Suppose, to the contrary, that $w$ is $M$-unmatched. If $w = u_2$, then by the maximality of $M$, the vertex $w$ is $M$-matched, a contradiction. Hence, $w \in V(T')$. If $vu_2$, then $M'$ has size~$|M| - 1 = \alpha'(T) - 1$. Since $\alpha'(T') = \alpha'(T) - 1$, the matching $M'$ is a maximum matching in $T'$. Applying the inductive hypothesis to the labeled tree $(T',S') \in \cT$ the vertex $w$ is $M'$-matched in $T'$ and therefore $M$-matched in $T$, a contradiction. Therefore, $vu_2 \notin M$, implying that $u_1u_2 \in M$ or $u_2u_3 \in M$. Renaming $u_1$ and $u_3$, if necessary, we may assume that $u_1u_2 \in M$. Thus, $M' = M \setminus \{u_1u_2\}$. Once again, $M'$ has size~$\alpha'(T) - 1$ and is a maximum matching in $T'$. Applying the inductive hypothesis to the labeled tree $(T',S') \in \cT$ the vertex $w$ is $M'$-matched in $T'$ and therefore $M$-matched in $T$, a contradiction. This establishes Property~(a).

We show secondly that Property~(b) holds. Let $w \in S_B \cup S_C$. Suppose that $w = u_1$. Let $M^*$ be a maximum matching  of $T'$. Then, the matching $M_w = M^* \cup \{u_2u_3\}$ is a maximum matching in $T$ such that $w$ is $M_w$-unmatched. Analogously, if  $w = u_2$, then there exists a maximum matching $M_w$ in $T$ such that $w$ is $M_w$-unmatched. Hence we may assume that $w \in V(T')$, for otherwise the desired result holds. Applying the inductive hypothesis to the labeled tree $(T',S') \in \cT$, there is a maximum matching $M^*$ of $T'$ such that the vertex $w$ is $M^*$-unmatched. Thus the matching $M_w = M^* \cup \{u_1u_2\}$ is a maximum matching in $T$ such that $w$ is $M_w$-unmatched. This establishes Property~(b).

We prove next that Property~(c) holds. As observed earlier, $\alpha'(T') = \alpha'(T) - 1$. In this case, we note that $\pc(T) = \pc(T') + 1$. Therefore, $\alpha'(T) + \frac{1}{2}\pc(T) = (\alpha'(T') + 1) + \frac{1}{2}(\pc(T') + 1) = \frac{n'}{2} + \frac{3}{2} = \frac{n}{2}$. This establishes Property~(c).

\medskip
\emph{Case~4. $(T,S)$ can be obtained from $(T',S')$ by operation $\cO_4$.} Let $v$ be the attacher in $T'$ and let $u_1u_2u_3u_4u_5$ be the path added to $T'$ and $vu_3$ the added edge. Then, $\sta(u_1) = \sta(u_5) = B$, $\sta(u_2) = \sta(u_4) = A$, and $\sta(u_3) = C$.

We show that $\alpha'(T') = \alpha'(T) - 2$. Since $M$ is a maximum matching in $T$, both vertices $u_2$ and $u_4$ are $M$-matched. If $vu_3 \notin M$, then $\alpha'(T') \ge |M'| = |M| - 2 = \alpha'(T) - 2$. If $vu_3 \in M$, then $M' = M \setminus \{u_1u_2,u_4u_5,vu_3\}$ and $|M'| = |M| - 3$. Applying the inductive hypothesis to the labeled tree $(T',S') \in \cT$, the matching $M'$ is not a maximum matching in $T'$ since the vertex $v \in S_A'$ is $M'$-unmatched. Hence in this case, $\alpha'(T') \ge |M'| + 1 = |M| - 2 = \alpha'(T) - 2$. In both cases, $\alpha'(T') \ge \alpha'(T) - 2$. Conversely, every matching in $T'$ can be extended to a matching in $T$ by adding to it the edges $u_1u_2$ and $u_4u_5$, and so $\alpha'(T) \ge \alpha'(T') + 2$. Consequently, $\alpha'(T') = \alpha'(T) - 2$.

We show firstly that Property~(a) holds. Let $w \in S_A$. Let $M$ be an arbitrary maximum matching in $T$ and let $M'$ be the restriction of $M$ to the tree $T'$. Suppose, to the contrary, that $w$ is $M$-unmatched. If $w = u_2$ or if $w = u_4$, then by the maximality of $M$, the vertex $w$ is $M$-matched, a contradiction. Hence, $w \in V(T')$. If $vu_3 \notin M$, then $M'$ has size~$|M| - 2 = \alpha'(T) - 2$. Since $\alpha'(T') = \alpha'(T) - 2$, the matching $M'$ is a maximum matching in $T'$. Applying the inductive hypothesis to the labeled tree $(T',S') \in \cT$, the vertex $w$ is $M'$-matched in $T'$  and therefore $M$-matched in $T$, a contradiction. Therefore, $vu_3 \in M$, implying that $M' = M \setminus \{u_1u_2,u_4u_5,vu_3\}$ and $|M'| = |M| - 3$.

We now consider the matching $M'$. Let $P_v$ be a longest path in $T'$ starting at $v$ and whose edges are alternately not in $M'$ and in $M'$ and whose vertices are alternately in $S_A'$ and in $S_C'$. Let $P_v$ be a $(v,z)$-path. By Observation~\ref{ob:cT}, every vertex in the set $S_A'$ has at least two neighbors in $S_B' \cup S_C'$. We now proceed analogously as in the proof of Case~1 presented earlier.

If $v = z$, then the vertex $v$ has a leaf-neighbor $u'$ that belongs to $S_B'$ and is $M$-unmatched. We now consider the matching $M^* = M' \cup \{u'v\}$. Since $M^*$ has size~$|M'| + 1 = |M| - 2 = \alpha'(T) - 2 = \alpha'(T')$, the matching $M^*$ is a maximum matching in $T'$. Since $w \notin \{u',v\}$ and $w$ is $M$-unmatched, the vertex $w$ is $M^*$-unmatched, contradicting Property~(a). Hence, $v \ne z$.

If $z \in S_C'$, then $P_v$ is an $M'$-augmenting $(v,z)$-path in $T'$. We now consider the matching $M^* = M \triangle E(P_v)$. Since $M^*$ has size~$|M'| + 1 = |M| - 2 = \alpha'(T) - 2 = \alpha'(T')$, the matching $M^*$ is a maximum matching in $T'$. Since $w \notin V(P_v)$ and $w$ is $M$-unmatched, the vertex $w$ is $M^*$-unmatched, contradicting Property~(a). Hence, $z \in S_A'$.

Since $z \in S_A'$, the vertex $z$ has a leaf-neighbor $u'$ that belongs to $S_B'$ and is $M$-unmatched. We now extend the path $P_v$ by adding to it the vertex $z'$ and the edge $zz'$ to produce a new path $P$ which is an $M'$-augmenting $(v,z')$-path in $T'$. We now consider the matching $M^* = M \triangle E(P)$. Since $M^*$ has size~$|M'| + 1 = |M| - 2 = \alpha'(T) - 2 = \alpha'(T')$, the matching $M^*$ is a maximum matching in $T'$. Since $w \notin V(P)$ and $w$ is $M$-unmatched, the vertex $w$ is $M^*$-unmatched, contradicting Property~(a). Therefore, $w$ is $M$-matched. Since $w$ was chosen to be an arbitrary vertex in $S_A$, this proves that Property~(a) holds.

We show secondly that Property~(b) holds. Let $w \in S_B \cup S_C$. Suppose that $w \notin V(T')$, and so $w \in \{u_1,u_3,u_5\}$. Let $M^*$ be a maximum matching of $T'$. If $w = u_1$, let $M_w = M^* \cup \{u_2u_3,u_4u_5\}$. If $w = u_3$, let $M_w = M^* \cup \{u_1u_2,u_4u_5\}$. If $w = u_5$, let $M_w = M^* \cup \{u_1u_2,u_3u_4\}$. In all cases, the matching $M_w$ is a maximum matching in $T$ such that $w$ is $M_w$-unmatched. Hence we may assume that $w \in V(T')$, for otherwise the desired result holds. Applying the inductive hypothesis to the labeled tree $(T',S') \in \cT$, there is a maximum matching $M^*$ of $T'$ such that the vertex $w$ is $M^*$-unmatched. Thus the matching $M_w = M^* \cup \{u_1u_2,u_4u_5\}$ is a maximum matching in $T$ such that $w$ is $M_w$-unmatched. This establishes Property~(b).

We prove next that Property~(c) holds. As observed earlier, $\alpha'(T') = \alpha'(T) - 2$. We note that here $\pc(T) = \pc(T') + 1$. Therefore, $\alpha'(T) + \frac{1}{2}\pc(T) = (\alpha'(T') + 2) + \frac{1}{2}(\pc(T') + 1) = \frac{n'}{2} + \frac{5}{2} = \frac{n}{2}$. This establishes Property~(c) and completes the proof of Lemma~\ref{lem:cT}.~\qed

\begin{lem}
If $T$ is a tree of order~$n \ge 3$ satisfying $\alpha'(T) + \frac{1}{2}\pc(T) = \frac{n}{2}$, then  $(T,S) \in \cT$ for some labeling $S$.
 \label{l:tree2}
\end{lem}
\proof We proceed by induction on the order $n \ge 3$ of a tree $T$ satisfying $\alpha'(T) + \frac{1}{2}\pc(T) = \frac{n}{2}$. If $n = 3$, then $T = P_3$. Letting $S = S^*_0$, we note that $(T,S) \in \cT$. This establishes the base case. Let $n \ge 4$ and assume that if $T'$ is a tree of order~$n'$, where $3 \le n' < n$, satisfying $\alpha'(T') + \frac{1}{2}\pc(T') = \frac{n'}{2}$, then $(T',S') \in \cT$ for some labeling $S'$. Let $T$ be a tree of order $n$ satisfying $\alpha'(T) + \frac{1}{2}\pc(T) = \frac{n}{2}$.

If $T$ is a star, then let $S$ be the labeling that assigns status $A$ to the central vertex of the star and status $B$ to every leaf. Then, $(T,S)$ can be obtained from the labeled tree $(P_3,S_0^*) \in \cT$ by repeated applications of Operation~$\cO_1$. Thus, $(T,S) \in \cT$.
Hence we may assume that $\diam(T) \ge 3$, for otherwise the desired result follows.

Let $P$ be a longest path in $T$ and suppose that $P$ is an $(r,u)$-path. Necessarily, $r$ and $u$ are leaves in $T$. We now root the tree $T$ at the vertex~$r$. Let $v$ be the parent of $u$, and let $w$ be the parent of $v$ in the rooted tree $T$. Since $\diam(T) \ge 3$, we note that $w \ne r$, implying that $d_T(w) \ge 2$. Let $x$ be the parent of~$w$.

Suppose that $d_T(v) \ge 4$. Let $T' = T - u$ and let $T'$ have order~$n'$, and so $n' = n - 1$. Then, $\alpha'(T) = \alpha'(T')$ and $\pc(T) = \pc(T') + 1$. Applying Lemma~\ref{l:tree1} to the tree $T'$, we see that
\[
\frac{n}{2} = \alpha'(T) + \frac{1}{2}\pc(T) = \alpha'(T') + \frac{1}{2}(\pc(T') + 1) \ge \frac{n'}{2} + \frac{1}{2} = \frac{n}{2}.
\]

Hence, we must have equality throughout the above inequality chain. In particular, $\alpha'(T') + \frac{1}{2}\pc(T') = \frac{n'}{2}$. Applying the inductive hypothesis to the tree $T'$, we note that $(T',S') \in \cT$ for some labeling $S'$. Let $S$ be the labeling obtained from $S'$ by assigning to the vertex~$u$ status~$B$. Since $v$ is a support vertex in $T'$, the vertex $v \in S'_A$ by Observation~\ref{ob:cT}(a). Therefore, the labeled tree $(T,S)$ can be obtained by the labeled tree $(T',S')$ by applying Operation~$\cO_1$, implying that $(T,S) \in \cT$. Therefore, we may assume that $d_T(v) \le 3$, for otherwise the desired result follows.

Suppose that $d_T(v) = 3$. Let $u_1$ and $u_2$ denote the two children of $v$, where $u = u_1$. We now let $T'$ be the tree obtained from $T$ by deleting $u_1,u_2$ and $v$; that is, $T' = T - D[v]$. Let $T'$ have order $n'$, and so $n' = n - 3$. Then, $\alpha'(T) = \alpha'(T') + 1$ and $\pc(T) = \pc(T') + 1$. Applying Lemma~\ref{l:tree1} to the tree $T'$, we see that
\[
\frac{n}{2} = \alpha'(T) + \frac{1}{2}\pc(T) = (\alpha'(T') + 1) + \frac{1}{2}(\pc(T')+1) \ge \frac{n'}{2} + \frac{3}{2} = \frac{n}{2}.
\]

Hence, we must have equality throughout the above inequality chain. In particular, $\alpha'(T') + \frac{1}{2}\pc(T') = \frac{n'}{2}$. Since $\{w,x\} \subseteq V(T')$, we note that $n' \ge 2$. If $n' = 2$, then the tree $T$ is determined (and is obtained from a star $K_{1,3}$ by subdividing one edge once). In this case, $n = 5$, $\alpha'(T) = 2$ and $\pc(T) = 2$, implying that $\alpha'(T) + \frac{1}{2}\pc(T) > \frac{n}{2}$, a contradiction. Therefore, $n' \ge 3$. Applying the inductive hypothesis to the tree $T'$, we note that $(T',S') \in \cT$ for some labeling $S'$. Let $S$ be the labeling obtained from $S'$ by assigning status~$A$ to the vertex~$v$ and status~$B$ to both $u_1$ and $u_2$. Further, if $w$ has status~$B$ in the labeling $S'$ (that is, if $w \in S_B'$), then change the status of $v$ from $B$ to $C$ in the labeling $S$ (that is, $w \in S_C$). Since the labeled tree $(T,S)$ can be obtained from the labeled tree $(T',S')$ by applying Operation~$\cO_3$, we therefore have that $(T,S) \in \cT$. Hence we may assume that $d_T(v) = 2$, for otherwise the desired result follows. Analogously, we may assume that every child of $w$ that is not a leaf has degree~$2$ in $T$.

Suppose that $d_T(w) \ge 3$. Let $w$ have $k$ children, and so $k = d_T(w) - 1 \ge 2$. By our earlier assumptions, every child of $w$ is a leaf or is a support vertex of degree~$2$ in $T$. Let $w$ have $k_1$ children that are non-leaves and $k_2$ children that are leaves. Thus, $k = k_1 + k_2$, $k \ge 2$, $k_1 \ge 1$, and $k_2 \ge 0$. Suppose that $k \ge 3$. Let $T' = T - \{u,v\}$ be obtained from $T$ by deleting the vertices $u$ and $v$, and let $T'$ have order~$n'$. Thus, $n' = n - 2$. Then, $\alpha'(T) = \alpha'(T') + 1$ and $\pc(T) = \pc(T') + 1$. Applying Lemma~\ref{l:tree1} to the tree $T'$, we see that
\[
\frac{n}{2} = \alpha'(T) + \frac{1}{2}\pc(T) = (\alpha'(T') + 1) + \frac{1}{2}(\pc(T') + 1) \ge \frac{n'}{2} + \frac{3}{2} > \frac{n}{2},
\]
a contradiction. Therefore, $k = 2$. Suppose that $k_1 = 1$, and so $k_2 = 1$. Let $v'$ be the leaf-neighbor of $w$. Let $T' = T - D[w]$ be obtained from $T$ by deleting the vertices $u, v, v'$ and $w$, and let $T'$ have order~$n'$. Thus, $n' = n - 4$. Then, $\alpha'(T) = \alpha'(T') + 2$ and $\pc(T) = \pc(T') + 1$. Applying Lemma~\ref{l:tree1} to the tree $T'$, we see that
\[
\frac{n}{2} = \alpha'(T) + \frac{1}{2}\pc(T) = (\alpha'(T') + 2) + \frac{1}{2}(\pc(T')+1) \ge \frac{n'}{2} + \frac{5}{2} > \frac{n}{2},
\]
a contradiction. Therefore, $k = k_1 = 2$. Let $v'$ be the child of $w$ different from $v$ and let $u'$ be the child of $v'$. We now let $T'$ be the tree obtained from $T$ by deleting $u,u',v,v'$ and $w$; that is, $T' = T - D[w]$. Let $T'$ have order $n'$, and so $n' = n - 5$. Then, $\alpha'(T) \ge \alpha'(T') + 2$ and $\pc(T) = \pc(T') + 1$. Applying Lemma~\ref{l:tree1} to the tree $T'$, we see that
\[
\frac{n}{2} = \alpha'(T) + \frac{1}{2}\pc(T) \ge (\alpha'(T') + 2) + \frac{1}{2}(\pc(T')+1) \ge \frac{n'}{2} + \frac{5}{2} = \frac{n}{2}.
\]

Hence, we must have equality throughout the above inequality chain. In particular, $\alpha'(T) = \alpha'(T') + 2$ and $\alpha'(T') + \frac{1}{2}\pc(T') = \frac{n'}{2}$. Since $P$ is a longest path in $T$, we note that $n' \ge 2$. If $n' = 2$, then the tree $T$ is determined (and is obtained from a star $K_{1,3}$ by subdividing every edge once). In this case, $n = 7$, $\alpha'(T) = 3$ and $\pc(T) = 2$, implying that $\alpha'(T) + \frac{1}{2}\pc(T) > \frac{n}{2}$, a contradiction. Therefore, $n' \ge 3$. Applying the inductive hypothesis to the tree $T'$, we note that $(T',S') \in \cT$ for some labeling $S'$. Recall that $x$ is the parent of $w$ in the rooted tree $T$. If $x \in S_B' \cup S_C'$, then by Lemma~\ref{lem:cT}, there exists a maximum matching $M'$ in $T'$ such that $w$ is $M'$-unmatched. In this case, $M' \cup \{xw,uv,u'v'\}$ is a matching in $T$ of size~$|M'| + 3 = \alpha'(T') + 3 > \alpha'(T)$, a contradiction. Therefore, $x \in S_A'$. Let $S$ be the labeling obtained from $S'$ by assigning status~$C$ to the vertex~$w$, status~$A$ to both $v$ and $v'$, and status~$B$ to both $u$ and $u'$. Since $x \in S_A'$, the labeled tree $(T,S)$ can be obtained by the labeled tree $(T',S')$ by applying Operation~$\cO_4$, implying that $(T,S) \in \cT$. Therefore, we may assume that $d_T(w) = 2$, for otherwise the desired result follows.

We now let $T'$ be the tree obtained from $T$ by deleting $u$ and $v$; that is, $T' = T - \{u,v\}$. Let $T'$ have order $n'$, and so $n' = n - 2$. Then, $\alpha'(T) = \alpha'(T') + 1$. Since $d_T(w) = 2$, we note that $w$ is a leaf in $T'$ and $\pc(T) = \pc(T')$. Applying Lemma~\ref{l:tree1} to the tree $T'$, we see that
\[
\frac{n}{2} = \alpha'(T) + \frac{1}{2}\pc(T) = (\alpha'(T') + 1) + \frac{1}{2}\pc(T') \ge \frac{n'}{2} + 1 = \frac{n}{2}.
\]

Hence, we must have equality throughout the above inequality chain. In particular, $\alpha'(T') + \frac{1}{2}\pc(T') = \frac{n'}{2}$.  If $n' = 2$, then the tree $T$ is determined and $T = P_4$. In this case, $n = 4$, $\alpha'(T) = 2$ and $\pc(T) = 1$, implying that $\alpha'(T) + \frac{1}{2}\pc(T) > \frac{n}{2}$, a contradiction. Therefore, $n' \ge 3$. Applying the inductive hypothesis to the tree $T'$, we note that $(T',S') \in \cT$ for some labeling $S'$. As observed earlier, $w$ is a leaf in $T'$. By Observation~\ref{ob:cT}(b), the vertex $w \in S_B'$. Let $S$ be the labeling obtained from $S'$ by assigning status~$A$ to the vertex~$v$, status~$B$ to the vertex~$u$ and reassigning the status of vertex $w$ to be status~$C$. Since $w \in S_B'$, the labeled tree $(T,S)$ can be obtained by the labeled tree $(T',S')$ by applying Operation~$\cO_2$, implying that $(T,S) \in \cT$. This completes the proof of Lemma~\ref{l:tree2_old}.~\qed

\medskip
As an immediate consequence of Lemmas~\ref{l:tree1},~\ref{lem:cT} and~\ref{l:tree2}, we have the result of Theorem~\ref{t:tree}. This result implies that given any graph $G$ of order~$n$, for every spanning tree $T$ of $G$ we know that $\alpha'(T) + \frac{1}{2}\pc(T) \ge \frac{n}{2}$. We are now in a position to prove the same lower bound holds for $G$ as well. Recall the statement of Theorem~\ref{t:thm2}.

\noindent \textbf{Theorem~\ref{t:thm2}} \emph{
If $G$ is a graph of order~$n$, then $\alpha'(G) + \frac{1}{2}\pc(G) \ge \frac{n}{2}$. Further for $n \ge 3$, $\alpha'(G) + \frac{1}{2}\pc(G) = \frac{n}{2}$ if and only if $G$ has a spanning tree $T$ such that \1
\\
\indent {\rm (a)} $(T,S) \in \cT$ for some labeling $S$.
\\
\indent {\rm (b)} $\alpha'(G) = \alpha'(T)$. \\
\indent {\rm (c)}  $\pc(G) = \pc(T)$.
}

\noindent \textbf{Proof.} By linearity, it suffices for us to restrict our attention to connected graphs. Let $G$ be a connected graph of order~$n$, and let $\cP$ be a minimum path cover in $G$, implying that $|\cP| = \pc(G)$. Let $T$ be a spanning tree obtained from the disjoint union of the $\pc(G)$ paths that belong to the path cover $\cP$ by adding $\pc(G) - 1$ edges (in such a way that the resulting graph is connected). Since $\cP$ is a path cover in $T$, we note that $\pc(T) \le |\cP| = \pc(G)$. Since adding edges to a tree cannot increase its path covering number, we note that $\pc(G) \le \pc(T)$. Consequently, $\pc(G) = \pc(T)$. Every matching in $T$ is a matching in $G$, implying that $\alpha'(G) \ge \alpha'(T)$. Hence, applying Lemma~\ref{l:tree1} to the tree $T$, we see that
\[
\alpha'(G) + \frac{1}{2}\pc(G) \ge \alpha'(T) + \frac{1}{2}\pc(T) \ge \frac{n}{2}.
\]

Suppose that $n \ge 3$ and $\alpha'(G) + \frac{1}{2}\pc(G) = \frac{n}{2}$. Then, we must have equality throughout the above inequality chain, implying that $\alpha'(G) = \alpha'(T)$ and $\alpha'(T) + \frac{1}{2}\pc(T) = \frac{n}{2}$. By Lemma~\ref{l:tree2}, $(T,S) \in \cT$ for some labeling $S$. Thus, conditions~(a),~(b) and~(c) in the statement of the theorem all hold. Conversely, if conditions~(a),~(b) and~(c) in the statement of the theorem all hold, then by Lemma~\ref{lem:cT}, $\alpha'(T) + \frac{1}{2}\pc(T) = \frac{n}{2}$, implying that
\[
\alpha'(G) + \frac{1}{2}\pc(G) = \alpha'(T) + \frac{1}{2}\pc(T) = \frac{n}{2}.
\]
This completes the proof of Theorem~\ref{t:thm2}.~\qed

\section{Applications to Domination Parameters}
\label{S:applic}


We remark that the upper bound of Theorem~\ref{Tdom_pc} on the total domination number in terms of its matching and path covering numbers is sharp even for graphs with minimum degree~$2$. For example, let $F$ be a star $K_{1,k-1}$ on $k \ge 3$ vertices and let $G$ be obtained from $F$ as follows: For each vertex $v$ of $F$, add a $6$-cycle and join $v$ to one vertex of this cycle. The resulting graph $G$ has order~$7k$, $\gamma_t(G) = 4k$, $\alpha'(G) = 3k + 1$ and $\pc(G) = k - 1$, implying that $\gamma_t(G) = 4k = \alpha'(G) + \pc(G)$.

As an application of Theorem~\ref{t:thm2}, we show, however, that the bound of Theorem~\ref{Tdom_pc} can be improved considerably if we restrict the minimum degree of the graph to be at least three.
If $G$ is a graph of order~$n$ with minimum degree at least~$3$, then it is well-known (see, for example,~\cite{Alfewy,HeYe_book}) that $\gamma_t(G) \le \frac{n}{2}$. Further the graphs achieving equality in this bound are characterized in~\cite{HeYe08}. We observe that these extremal graphs in~\cite{HeYe08} do not achieve equality in the bound of Theorem~\ref{t:thm2}. Therefore as a consequence of these results on total domination in graphs with minimum degree at least~$3$, the following improvement of Theorem~\ref{Tdom_pc} follows readily from Theorem~\ref{t:thm2}: If $G$ is a graph with $\delta(G) \ge 3$, then $\gamma_t(G) < \alpha'(G) + \frac{1}{2}\pc(G)$. We remark that there are two (infinite) families of connected cubic graphs $G$ satisfying $\gamma_t(G) = \frac{n}{2}$ (as shown in~\cite{HeYe08}). Each such graph $G$ satisfies $\alpha'(G) + \frac{1}{2}(\pc(G) - 1) = \frac{n}{2}$. Consequently, there are infinitely many connected cubic graphs $G$ satisfying $\gamma_t(G) = \alpha'(G) + \frac{1}{2}(\pc(G) - 1)$. This establishes the result of Corollary~\ref{t:thm4}. Recall its statement.

\medskip
\noindent \textbf{Corollary~\ref{t:thm4}} \emph{
If $G$ is a graph with $\delta(G) \ge 3$, then \[
\gamma_t(G) \le \alpha'(G) + \frac{1}{2}(\pc(G) - 1),\]
and this bound is tight.
}

As observed earlier, the relationship between the neighborhood total domination number and the matching number of a graph behaves quite differently from the relationship between the domination number and the matching number of a graph. Certainly there exist many graphs $G$ for which $\alpha'(G) \ge \gnt(G)$. For example, if $G$ is obtained from a connected graph $H$ by adding at least one pendant edges to each vertex of $H$, then $\alpha'(G) = \gnt(G) = n(H)$ (noting that the set $V(H)$ is a minimum NTD-set in $G$). If $G = K_{n,n}$ where $n \ge 3$, then $\alpha'(G) = n > 2 = \gnt(G)$. What is not quite as obvious is that there exist an infinite number of graphs for which $\alpha'(G) < \gnt(G)$, even for arbitrary large (but fixed) minimum degree.

\begin{thm}
For every integer $\delta \ge 1$, there exists a graph $G$ with $\delta(G) = \delta$ satisfying $\gnt(G) > \alpha'(G)$.
\label{t:thm1}
\end{thm}
\proof When $\delta = 1$, take $G$ to be a star on at least two vertices, while if $\delta = 2$, take $G$ to be a $5$-cycle. Hence we may assume that $\delta \ge 3$. Let $n = \delta(\delta + 1)$, let $A$ be a set of $n$ elements, and let $B$ be the set of all $\delta$-element subsets of $A$. Viewing the elements of $A$ and $B$ as vertices, let $G = G_n^{\delta}$ be the bipartite graph formed by taking $A$ as one partite set and $B$ as the other partite set where vertex $v_a$ corresponding to element $a\in A$ is adjacent to vertex $v_b$ corresponding to the $\delta$-element subset $b\in B$ if and only if $a$ is contained in $b$. Thus, every vertex in $B$ has degree~$\delta$, while every vertex in $A$ has degree~${n-1 \choose \delta-1}$. Furthermore, $G$ is a bipartite graph with minimum degree~$\delta$ and order $n + {n \choose \delta}$. Now, $\alpha'(G) \le \min(|A|,|B|) = |A| = n$. It is easy to find a matching of size~$n$ in $G$, and so $\alpha'(G) = n$.

We show next that $\gnt(G) \ge n+1$. Let $D$ be an arbitrary NTD-set in $G$ and write $D_A = D \cap A$ and $D_B = D \cap B$. It follows that $|D| = |D_A| + |D_B|$.
We show that $|D_B| \ge \delta + 1$. Let $S$ be an arbitrary $\delta$-element subset of $A$ and let $v_S$ be the vertex in $B$ associated with $S$. Suppose that $D_B$ contains no neighbor of a vertex in $S$; that is, $D_B \cap N(S) = \emptyset$. On the one hand, if $S$ contains a vertex of $A \setminus D_A$, then such a vertex would not be dominated by $D$. On the other hand, if $S \subseteq D_A$, then the vertex $v_S$ would be isolated in the subgraph induced by $N(D)$. Both cases produce a contradiction. Therefore, $D_B$ contains at least one vertex in $N(S)$. This is true for any $\delta$-element subset, $S$, of $A$. In particular, this implies that $|D_B| \ge \lfloor |A| / \delta \, \rfloor = \lfloor n / \delta \, \rfloor = \delta + 1$ since each vertex of $B$ dominates exactly~$\delta$ vertices in $A$ and $|A| = \delta^2 + \delta$.

We are now in a position to show that $|D| \ge n+1$. If $|D_A| = n - k$ for some $k$ where $k < \delta$, then $|D| = |D_A| + |D_B| > (n - \delta) + (\delta + 1) = n+1$. Hence, we may assume that $|D_A| = n - k$ for some $k$ where $k \ge \delta$, for otherwise $|D| \ge n+1$, as desired.  In order to dominate the ${k \choose \delta}$ vertices in $B$ associated with $\delta$-element subsets of $A \setminus D_A$, all these ${k \choose \delta}$ vertices belong to the set $D_B$. Further, since $D_B$ contains at least one vertex in $N(S)$ for every $\delta$-element subset $S$ of $A$, the set $D_B$ contains at least $\lfloor |D_A| / \delta \, \rfloor = \lfloor (n-k) / \delta \, \rfloor$ additional vertices, implying that
\begin{equation}
|D| \ge n - k + {k \choose \delta} + \left\lfloor \frac{n-k}{\delta} \, \right\rfloor.
\label{Eq1}
\end{equation}

Let $f(k)$ denote the function on the right hand side of Inequality~(\ref{Eq1}). It suffices for us to show that $f(k) \ge n+1$. If $k = \delta$, then $f(k) = n+1$. Hence we may assume that $k \ge \delta + 1$. Recall that $n = \delta(\delta + 1)$ and, by assumption, $\delta \ge 3$. If $D_A = \emptyset$, then $k = n = \delta(\delta + 1)$ and $f(k) = {k \choose \delta} = {\delta (\delta + 1) \choose \delta} > n+1$. Suppose that $D_A \ne \emptyset$, and so $n - k \ge 1$. If $n - k < \delta$, then $k \ge \delta^2 + 1$ and $f(k) \ge 1 + {\delta^2 + 1 \choose \delta} \ge n + 1$. If $n - k \ge \delta$  and $k \ge \delta + 1$, then $\left\lfloor \frac{n-k}{\delta} \, \right\rfloor \ge 1$ and ${k \choose \delta} \ge k$, implying that $f(k) \ge (n-k) + k + 1 = n + 1$. In all cases, we have shown that $f(k) \ge n+1$, as claimed. Therefore, if $k \ge \delta$, then $|D| \ge f(k) \ge n+1$.
Since $D$ is an arbitrary NTD-set of $G$ and $|D| \ge n+1$, we deduce that $\gnt(G) \ge n+1 > \alpha'(G)$, as desired.~\qed

\medskip
Although we have seen that, given any graph $G$, $\gnt(G)$ may be larger than $\alpha'(G)$, it turns out that we can upper bound $\gnt(G)$ in terms of the matching number and the path covering number of $G$.  As a consequence of Theorem~\ref{t:oddn} and Theorem~\ref{t:thm2}, we have a proof of Corollary~\ref{t:thm3}. Recall its statement.

\medskip
\noindent \textbf{Corollary~\ref{t:thm3}} \emph{
If $G$ is a connected graph on at least three vertices, then
\[
\gnt(G) \le \alpha'(G) + \frac{1}{2}\pc(G)
\]
unless $G \in \{P_3,P_5,C_5\}$ in which case $\gnt(G) = \alpha'(G) + \frac{1}{2}(\pc(G)+1)$.
}

\noindent \textbf{Proof.} Let $G$ be a connected graph of order $n \ge 3$. If $G = C_5$ or if $G$ is a subdivided star, then by Theorem~\ref{t:oddn}, $\gnt(G) = (n+1)/2$. If $G$ is a subdivided star and $n \ge 7$, then $\alpha'(G) = (n-1)/2$ and $\pc(G) \ge 2$, and so $\alpha'(G) + \frac{1}{2}\pc(G) \ge (n+1)/2 = \gnt(G)$. If $G$ is a subdivided star and $n < 7$, then $G \in \{P_3,P_5\}$ and $\alpha'(G) + \frac{1}{2}(\pc(G)+1) \ge (n+1)/2 = \gnt(G)$. If $G = C_5$, then once again $\alpha'(G) + \frac{1}{2}(\pc(G)+1) \ge (n+1)/2 = \gnt(G)$. Therefore, we may assume that $G$ is neither a $5$-cycle nor a subdivided star, for otherwise the desired result follows. Hence by Theorem~\ref{t:oddn}, $\gnt(G) \le \frac{n}{2}$ and the desired result follows by Theorem~\ref{t:thm2}.~\qed

\medskip
We remark that the bound of Corollary~\ref{t:thm3} is tight as illustrated by the following example. Let $H$ be an arbitrary hamiltonian graph of odd order~$k \ge 1$ and let $G$ be the graph obtained from $H$ by the following operation: for each vertex $x \in V(H)$, add a new path $P_3$ and an edge joining its central vertex to $x$. Equivalently, for $k \ge 1$ odd, the graph $G$ is obtained from the disjoint union of $k \ge 1$ stars $K_{1,3}$ by selecting one leaf from each star and adding any number of edges between the selected $k$ leaves so that they induce a graph that contains a hamiltonian path. The resulting graph $G$ has order~$4k$, $\gnt(G) = 2k$, $\alpha'(G) = k + (k-1)/2$ and $\pc(G) = k + 1$, implying that $\gnt(G) = 2k = \alpha'(G) + \frac{1}{2}\pc(G)$.


\medskip
\begin{thebibliography}{99}

\bibitem{Alfewy} D. Archdeacon, J. Ellis-Monaghan, D. Fischer, D. Froncek, P.C.B. Lam, S. Seager, B. Wei, and R. Yuster, Some remarks on domination. \textit{J. Graph Theory} \textbf{46} (2004), 207--210.



\bibitem{ArSi11} S. Arumugam, C. Sivagnanam, Neighborhood total
    domination in graphs. \textit{Opuscula Mathematica} \textbf{31} (2011), 519--531.

\bibitem{BoCo79} Bollob\'{a}s B. and  E. J. Cockayne, Graph-theoretic parameters concerning domination, independence, and irredundance. \textit{J. Graph Theory} \textbf{3} (1979), 241--249.


\bibitem{BrCaVi00} R. C. Brigham, J. R. Carrington, and R. P.
    Vitray, Connected graphs with maximum total domination number.
    \textit{J. Combin. Comput. Combin. Math.} \textbf{34} (2000), 81--96.

\bibitem{ChMc} V. Chv\'{a}tal and C. McDiarmid, Small transversals in hypergraphs. \textit{Combinatorica} \textbf{12} (1992), 19--26.


\bibitem{CoDaHe80} E. J. Cockayne, R. M. Dawes, and S. T. Hedetniemi, Total domination in graphs. \textit{Networks} \textbf{10} (1980),
    211--219.

\bibitem{CoHeSl79} E. J. Cockayne, S. T. Hedetniemi and P. J. Slater, Matchings and transversals in hypergraphs, domination and independence-in trees. \textit{J. Combin. Theory B} \textbf{27} (1979), 78--80.

\bibitem{DeLa07}  E. DeLaVi\~{n}a, Q. Liu, R. Pepper, B. Waller, and D. B. West, Some conjectures of Graffiti.pc on total domination. \textit{Congressus Numer.} \textbf{185} (2007), 81--95.

\bibitem{DoHaHeMa00} G. Domke, J. H. Hattingh, M. A. Henning, and L. Markus, Restrained domination in graphs with minimum degree two. \textit{J. Combin. Math. Combin. Comput.} \textbf{35} (2000), 239--254.




\bibitem{hhs1} T. W. Haynes, S. T. Hedetniemi, and P. J. Slater,
    \emph{Fundamentals of Domination in Graphs}, Marcel Dekker, Inc.
    New York, 1998.

\bibitem{hhs2} T. W. Haynes, S. T. Hedetniemi, and P. J. Slater
    (eds), \emph{Domination  in Graphs: Advanced Topics}, Marcel
    Dekker, Inc. New York, 1998.

\bibitem{He09} M. A. Henning, Recent results on total domination in graphs: A survey. \textit{Discrete Math}. \textbf{309} (2009),
    32--63.


\bibitem{HeKaShYe08} M. A. Henning, L. Kang, E. Shan  and A. Yeo, On matching and total domination in graphs. \textit{Discrete Math.} \textbf{308} (2008), 2313–-2318.


\bibitem{HeRa13} M. A. Henning and N. J. Rad, Bounds on neighborhood total domination in graphs. \textit{Discrete
    Applied Math.} \textbf{161} (2013), 2460--2466.

\bibitem{HeWa14} M. A. Henning and K. Wash, Trees with large neighborhood total domination number, mansucript (http://arxiv.org/abs/1408.0109).


\bibitem{HeYe06} M. A. Henning and A. Yeo, Total domination and matching numbers in claw-free graphs. \textit{Electronic J. Combin.} \textbf{13} (2006), \#59.

\bibitem{HeYe08} M. A. Henning and A. Yeo, Hypergraphs with large transversal number and with edge sizes
            at least three. \textit{J. Graph Theory}. \textbf{59} (2008), 326--348.

\bibitem{HeYe11} M. A. Henning and A. Yeo, Perfect matchings in total domination critical graphs. \textit{Graphs Combin.}
    \textbf{27} (2011), 685--701.

\bibitem{HeYe13} M. A. Henning and A. Yeo, Total domination and matching numbers in graphs with all vertices in triangles. \textit{Discrete Math.} \textbf{313} (2013), 174–-181.

\bibitem{HeYe_book} M. A. Henning and A. Yeo, \emph{Total domination in graphs (Springer Monographs in Mathematics)} 2013. ISBN: 978-1-4614-6524-9 (Print) 978-1-4614-6525-6 (Online).

\bibitem{HeYe_chapter}  M. A. Henning and A. Yeo, Total domination number versus matching number. \emph{Total domination in graphs (Springer Monographs in Mathematics)} (2013), 77--81.


\bibitem{OWest10} S. O and D. B. West, Balloons, cut-edges, matchings, and total domination in regular graphs of odd degree. \textit{J. Graph Theory} \textbf{64} (2010), 116–-131.

\bibitem{Pl03} M. Plummer, Factors and Factorization. 403--430. \emph{Handbook of Graph Theory} \emph{ed.} J. L. Gross and J. Yellen. CRC Press, 2003, ISBN: 1-58488-092-2.

\bibitem{Pu95} W. R. Pulleyblank, Matchings and Extension. 179--232. \emph{Handbook of Combinatorics} \emph{ed.} R. L. Graham, M. Gr\"{o}tschel, L. Lov\'{a}sz. Elsevier Science B.V. 1995, ISBN  0-444-82346-8.

\bibitem{ShChCh10} W. C. Shiu, X. Chen, W. H. Chan, Some results on matching and total domination in graphs. \textit{Appl. Anal. Discrete Math.} \textbf{4} (2010), 241--252.

\bibitem{WaKaSh09} H. Wang, L. Kang, and E. Shan, Matching properties on total domination vertex critical graphs. \textit{Graphs Combin.} \textbf{25}(6) (2009), 851--861.

\end{thebibliography}
\end{document}